\documentclass{article}

\usepackage{amsmath}
\usepackage{amsthm}
\usepackage{graphicx}
\usepackage{upgreek}
\usepackage{mathrsfs}
\usepackage{xcolor}
\usepackage{bm,bbm}
\usepackage{amsfonts}
\usepackage{amssymb}
\usepackage{csquotes}
\usepackage{stmaryrd}
\usepackage{ulem}
\usepackage{caption}
\usepackage{subcaption}

\usepackage[
pdftitle={HJB G-VMC},
pdfcreator={pdflatex},
pdfsubject={preprint},
hyperindex = {true},
colorlinks = {true},
linkcolor = {blue},
citecolor = {blue}
]{hyperref}

\usepackage{enumerate}

\usepackage{soul}
\usepackage[margin=1.1in]{geometry}
\numberwithin{equation}{section}

\newtheorem{proposition}{Proposition}[section]

\newtheorem{example}{Example}[section]

\newtheorem{remark}{Remark}[section]

\newcommand{\EE}{{\mathbb E}}
\newcommand{\FF}{{\mathbb F}}
\newcommand{\PP}{{\mathbb P}}
\newcommand{\R}{\mathbb{R}}

\newcommand{\calC}{{\mathcal C}}

\newcommand{\calF}{{\mathcal F}}

\newcommand{\sig}{\sigma}

\newcommand{\kap}{\kappa}

\newcommand{\diag}{{\rm diag}}

\newcommand{\noi}{\noindent}

\newcommand{\footremember}[2]{%
    \footnote{#2}
    \newcounter{#1}
    \setcounter{#1}{\value{footnote}}%
}
\newcommand{\footrecall}[1]{%
    \footnotemark[\value{#1}]%
} 

\title{
Quantum-inspired nonlinear Galerkin ansatz for high-dimensional HJB equations
}
\author{%
  Chuhao Sun\footremember{UM}{Department of Mathematics, University of Michigan, Ann Arbor, MI 48109}%
  \and Asaf Cohen\footrecall{UM}%
  \and James Stokes\footrecall{UM}
  \and Shravan Veerapaneni\footrecall{UM} 
  }
\begin{document}
\maketitle
\begin{abstract}
Neural networks are increasingly recognized as a powerful numerical solution technique for partial differential equations (PDEs) arising in diverse scientific computing domains, including quantum many-body physics. In the context of time-dependent PDEs, the dominant paradigm involves casting the approximate solution in terms of stochastic minimization of an objective function given by the norm of the PDE residual, viewed as a function of the neural network parameters. Recently, advancements have been made in the direction of an alternative approach which shares aspects of nonlinearly parametrized Galerkin methods and variational quantum Monte Carlo, especially for high-dimensional, time-dependent PDEs that extend beyond the usual scope of quantum physics. This paper is inspired by the potential of solving Hamilton-Jacobi-Bellman (HJB) PDEs using Neural Galerkin methods and commences the exploration of nonlinearly parametrized trial functions for which the evolution equations are analytically tractable. As a precursor to the Neural Galerkin scheme, we present trial functions with evolution equations that admit closed-form solutions, focusing on time-dependent HJB equations relevant to finance.
\end{abstract}

\noi{\bf Keywords:}   Hamilton-Jacobi-Bellman equations, utility indifference pricing

\noi{\bf AMS Classification:}

\tableofcontents
\section{Introduction}
Neural networks are playing an increasingly important role in the solution of scientific and engineering problems that have resisted attack by traditional numerical methods. One arena that has received considerable attention in recent years is the solution of partial differential equations (PDEs). The central paradigm revolves around physics-informed neural networks (PINNs) \cite{raissi2019physics} and the deep Galerkin methods (DGM) \cite{sirignano2018dgm}, whose key insight is to treat the norm of the PDE residual as a stochastic objective function for optimizing the weights of a deep neural network. Recently an alternative approach has emerged \cite{reh2022variational, bruna2022neural, zhao2022quantum}, which is suitable specifically for time-dependent PDEs. Because this approach  is closely related to both Galerkin methods in traditional numerical analysis and also shares ideas from time-dependent variational quantum Monte Carlo with neural-network quantum states \cite{carleo2017solving}, it has been referred to as Neural Galerkin \cite{bruna2022neural} and alternatively as VMC for PDEs \cite{reh2022variational, zhao2022quantum}. Since the essential ingredient for the purposes of this work is the nonlinear parametrization of the trial function, this work adopts the acronym NG to refer to Nonlinear and Neural Galerkin interchangeably (the latter being a special case of the former). A pervasive challenge shared across PINNs, DGM, and NG is the lack of interpretability of the neural network architectures representing the PDE solution, making the task of architecture selection for specific PDEs rather arbitrary. In an effort to alleviate the interpretability problem, and aid the development of solid design principles, this work focuses on aspects of the NG method which are amenable to analytical control. In particular, we develop trial functions for which the evolution equations can be solved in closed form, focusing on time-dependent Hamilton-Jacobi-Bellman (HJB) equations arising in finance. Specifically, we identify a simple class of interpretable trial fucntions and apply them to the HJB equation in diverse dimensions. In addition to providing an approximate closed-form solution method of independent interest, the proposed trial functions can be used as the starting point of neural-network architectures based on normalizing flows \cite{papamakarios2021normalizing}. The overarching strategy outlined here for HJB equations departs from the prevalent one-size-fits-all architecture approach and holds potential for wider applications in designing neural network architectures tailored for specific PDE problems.

The financial framework we select for illustrating our methodology revolves around utility maximization and utility indifference pricing of options in incomplete markets. In complete markets, derivative pricing relies on replication and hedging strategies. However, in incomplete markets where perfect replication is unfeasible, a unique price may not exist. Consequently, utility indifference pricing comes into play, representing the price at which an investor attains an equivalent expected utility level, whether they purchase the derivative or not. We examine scenarios where the derivative's value is tied to non-tradable assets, allowing investors to allocate their initial capital strategically, depending on their choice to buy at the indifferent price. This approach necessitates solving two utility maximization problems.

This paper is structured as follows. In Section \ref{sec:2} we review the NG method, paying specific attention to its quantum-mechanical aspects. Section \ref{sec:3} provides a short introduction to HJB equations and their applications in mathematical finance. In Section \ref{sec:4} we apply the NG method to the financial models. Finally, section \ref{sec:conclusion} presents conclusions and future directions.

\section{The Nonlinear/Neural Galerkin Method}\label{sec:2}

In this section we abstract and synthesize some of the results from \cite{reh2022variational, bruna2022neural, zhao2022quantum}. Here, we will use $\mathbb{F}$ to denote either $\mathbb{R}$ or $\mathbb{C}$. Given a subset $\Omega \subseteq \mathbb{R}^d$, $L^2(\mathbb{R}^d)$ denotes the space of square-integrable functions with inner product $\langle \cdot | \cdot \rangle$ defined with the following conventions $\langle f_1 | f_2 \rangle = \int_{\Omega} {\rm d}^dx \overline{f_1(x)} f_2(x)$ and $\Vert \cdot \Vert$ denotes the induced norm.

The NG method targets a time-dependent PDE initial value problem of the following general form,
\begin{equation}\label{e:IVP}
    \left\{\begin{aligned}
    (\partial_t u)(t,x) & = F\big[t,x,u(t,x),(D^{\alpha_1} u)(t,x),\ldots, (D^{\alpha_m} u)(t,x)\big], & (t,x) \in [0,T] \times \Omega, \\
    u(0,x) & = f(x), & x \in \Omega.
    \end{aligned}\right.
\end{equation}
where $\Omega \subseteq \mathbb{R}^d$ and $F : \mathbb{R}^{d+1}\times\mathbb{F}^{\kap+1} \longrightarrow \mathbb{F}$ is an arbitrary function of the space-time coordinates, the state variable, and $\kap \geq 1$ spatial derivatives of arbitrary order $\alpha_1,\ldots,\alpha_\kap \in \mathbb{N}^d$, where the mixed partial derivative of order $\alpha = (\alpha^1,\ldots,\alpha^d) \in \mathbb{N}^d$ is defined as
\begin{align*}
    (D^{\alpha}u)(t,x) & := \frac{\partial^{|\alpha|}}{\partial x_1^{\alpha^1}\cdots \partial x_d^{\alpha^d}}u(t,x_1,\ldots,x_d), \\
    |\alpha| & := \alpha^1 + \cdots + \alpha^d.
\end{align*}

In addition to a time-dependent PDE, the NG assumes as input a family of sufficiently regular functions $u_\theta : \Omega \longrightarrow \mathbb{F}$, which have the property that they satisfy the boundary conditions on $\partial\Omega$ and are differentiably parametrized by $\theta \in \mathbb{R}^{p+1}$. The initial parameters $\theta(0)\in \mathbb{R}^{p+1}$ are chosen such that $u_{\theta(0)}$ approximates the initial condition within the desired tolerance $u_{\theta(0)}\cong f$ and then the subsequent values $\theta(t) \in \mathbb{R}^{p+1}$ are determined for $t>0$ by solving
\begin{equation}\label{e:general}
    M(\theta(t)) \, \theta'(t) = V(t,\theta(t)),
\end{equation}
where for $1 \leq i,j \leq p+1$,
\begin{equation*}
    M_{ij}(\theta) 
    := \operatorname{Re}\left[\left\langle \frac{\partial u_\theta}{\partial \theta_i} \middle| \frac{\partial u_\theta}{\partial \theta_j} \right\rangle\right], \quad \quad 
    V_i(t,\theta) := \operatorname{Re}\left[\left\langle \frac{\partial u_\theta}{\partial \theta_i} \middle| \widetilde F(t,\theta,\cdot) \right\rangle\right],
\end{equation*}
and
\begin{equation*}
    \widetilde F(t,\theta,x) := F\big[t,x,u_\theta(x),(D^{\alpha_1} u_\theta)(x),\ldots, (D^{\alpha_n} u_\theta)(x)\big].
\end{equation*}

In order for the integrals defining $M$ and $V$ to converge, by Cauchy-Schwarz inequality it suffices to require that $\widetilde{F}(t,\theta,\cdot) \in L^2(\Omega)$ and $\frac{\partial u_\theta}{\partial \theta_i}\in L^2(\Omega)$ for all $t \in [0,T]$, $\theta\in\mathbb{R}^{p+1}$ and $1\leq i \leq p+1$. 

It is convenient (but not necessary) to also demand that $u_\theta \in L^2(\Omega)$ for all $\theta \in \mathbb{R}^{p+1}$. If this condition holds, then without loss of generality we can express $u_\theta$ in terms of a unit vector $\psi_{\bm{\theta}} \in L^2(\Omega)$ parametrized by $\bm{\theta} \in \mathbb{R}^{p}$ as follows,
\begin{equation*}
    u_\theta(x) = \exp(\theta_0) \psi_{\bm{\theta}}(x) , \quad \quad \quad x \in \Omega,
\end{equation*}
where now $\theta = (\theta_0,\bm{\theta}) \in \mathbb{R}^{p+1}$ and the first component simply controls the norm of $u_\theta$, since
\begin{equation*}
    \Vert u_\theta \Vert = \exp(\theta_0) \Vert \psi_\theta \Vert = \exp(\theta_0).
\end{equation*}

In order to render the NG scheme practical for high-dimensional trial functions $u_\theta$ described by neural networks, numerical estimation of the high-dimensional integrals defining $M$ and $V$ is necessary. Proposals based on active learning \cite{bruna2022neural} and normalizing flows \cite{reh2022variational} have been put forth. In this work, we bypass the sampling problem by focusing on trial functions for which the HJB integrals can be computed in closed form.

\section{Applications of HJBs in mathematical finance: utility indifference pricing and utility optimization problems}\label{sec:3}

In this section, we delve into HJB equations and their role in mathematical finance. We commence with a concise overview of HJB equations (Section \ref{sec:HJB}), and subsequently, explore their applications in utility indifference pricing of derivatives within incomplete markets featuring exponential utility functions, along with utility optimization problems. The concept of utility indifference pricing is elaborated further in Section \ref{sec:indifference}. We then introduce the underlying stochastic model in Section \ref{sec:model}. Moving forward, Section \ref{sec:uti} presents an introduction to utility indifference pricing and utility optimization. Section \ref{sec:hjb} supplies the HJB equation in its general form. Lastly, in Section \ref{sec:eg}, we furnish a few illustrative examples where HJB equations exhibit simplified expressions.

\subsection{HJB equations}\label{sec:HJB}
The HJB equation is a mathematical expression that serves as the optimality criterion for a range of control problems. It originates from the dynamic programming principle and finds application in diverse domains, encompassing classical deterministic control problems as well as stochastic control problems. Generally, consider a control problem on the time interval $[0,T]$ with the following reward function given at time $t \in [0,T]$,
\begin{align*}
    J(t,x,u) = \EE\Big[\int_t^T h(s,X_s^u,u_s)ds + g(X^u_T)\;\Big|\;X^u_t = x\Big],
\end{align*}
with the diffusion process $X^u = (X^u_s)_{s\in[0,T]}$ given by
\begin{align*}
    X^u_t = X_0^u+\int_0^tb(X^u_s,u_s)ds+ \int_0^t\sig(X^u_s,u_s)dW_s, \quad t \in [0,T].
\end{align*}
Here, $T >0$ is a fixed time horizon, the function $h$ is called the running reward and $g$ is the terminal reward, $u = (u_t)_{t\geq 0}$ is a control process, taking values in some space, and is chosen by a decision maker to manipulate the process $X^u$, $(W_s)_{s\in[0,T]}$ is a standard Brownian motion. The function $J(t,x,u)$ is the expected return, given the initial time $t$, initial state $X^u_t = x$, and the control $u$. The control is assumed to satisfy some technical constraints and is called an {\it admissible control}.

Set the {\it value function} at time $t$, when the initial state is $X^u_0=x$, by:
\begin{align*}
    V(t,x) = \sup_u J(t,x,u).
\end{align*}
This is the optimal reward over all admissible controls $u$. Under some conditions on the data of the problem, the value function $V$ is the unique solution of the following PDE, called the HJB equation, \cite{Huyen}
\begin{align*}
\begin{cases}
    (\partial_t\phi)(t,x) + H(t,x,(D_x\phi)(t,x),(D^2_{x}\phi)(t,x)) = 0,\qquad (t,x) \in [0,T]\times\R^d,\\
    \phi(T,x) = g(x),
    \end{cases}
\end{align*}
where $D_x\phi$ is the gradient of $\phi$ and $D^2_x\phi$ is the Hessian matrix;
the function $H$ is the {\it Hamiltonian} associated with the problem, where in the case described above, is given by:
\begin{align*}
    H(t,x,p,m) = \sup_u\Big\{b(x,u)\cdot p +\frac 12 \text{tr}(\sig(x,u)\sig^T(x,u) m) + h(x,t,u)\Big\}.
\end{align*}

\subsection{Utility indifference pricing}\label{sec:indifference}

We choose to illustrate our methods by focusing on HJB equations within the context of financial models, specifically addressing utility maximization problems and utility indifference pricing of options in incomplete markets. In a complete market, derivative pricing is achieved through replication and hedging. However, within an incomplete market, where perfect replication is unattainable, a unique price may not exist. Thus, utility indifference pricing comes into play. This represents the current derivative price at which an investor reaches the same expected utility level, whether they purchase the derivative or not.

In cases where the derivative is reliant on assets that cannot be traded, the investor retains the flexibility to utilize their remaining initial capital, contingent on their decision to buy the derivative at the indifferent price (considering a buy price at the moment). This allows them to engage in trading with other tradable assets to maximize their terminal utility. To address this complex scenario, one must tackle two utility maximization problems.

We provide the HJB equation for a general optimal utility function of the investor and test our method in several examples. The reader is referred to \cite{survey} for a survey on utility indifference pricing.

\subsection{A generic stochastic model in mathematical finance}\label{sec:model}
Fix a time horizon $T > 0$ and consider a filtered probability space $(\Omega, \FF,\PP,(\calF_t)_{0 \leq t\leq T})$ supporting $B_t=(B^1_t,\dots,B^m_t)$ and $W_t = (W^1_t,\dots,W^n_t)$, $m$-dimensional and $n$-dimensional standard Brownian motions, respectively, with correlation $<B^i_t,W^j_t> = \rho^{i,j}\in[-1,1]$, $i\in\{1,\ldots,m\}$ and $j\in\{1,\ldots,n\}$. 
The probability space also supports an $m$-dimensional processes $S = (S_t)_{t\in[0,T]}$, which represents the tradable risky assets, and an $n$-dimensional non-tradable process $Y = (Y_t)_{t\in[0,T]}$ (to which we also refer to as the {\it auxiliary process}) satisfying the following stochastic differential equation (SDE) on $[0,T]$:
\begin{align}
\notag
    dS_t &= \mu(t,S_t,Y_t)dt+\sigma(t,Y_t,S_t)dB_t,\quad S_0 = s_0>0,\\
    dY_t &= b(t,S_t,Y_t)dt + a(t,S_t,Y_t)dW_t,\quad Y_0 = y_0 > 0, \label{dyn-Y}
\end{align}
where the functions $\mu:[0,T]\times\R^m\times\R^n\to\R^m$, $\sigma:[0,T]\times\R^m\times\R^n\to\R^{m\times m}$, $b:[0,T]\times\R^m\times\R^n\to\R^n$, and $a:[0,T]\times\R^m\times\R^n\to\R^{n\times n}$ guarantee a unique strong solution to the $(m+n)$-dimension SDE above. 

Consider also a risk-free asset $R$, whose value at time $t$ is $R_t = R_0e^{rt}$, where $r$ is the constant interest rate and $R_0 >0$ is a constant.

Let $\theta = (\theta_t)_{t\in[0,T]}$ denote the portfolio vector the investor holds at time $t$, i.e., $\theta_t^i$ is the amount of asset $S^i$ the investor holds. Then, the wealth process with an initial wealth $x_0$ satisfies on $[0,T]$:
\begin{align}\label{wealth}
    dX^{\theta,x_0}_t = \theta_tdS_t + r_t(X^{\theta,x_0}_t-\theta_tS_t)dt, \qquad X^{\theta,x_0}_0 = x_0,
\end{align}
where the first term on the right-hand side denotes the change in the risky investment and the second term denotes the change in the risk-free investment. We call $\theta$ an admissible control if it is progressively measurable with respect to (w.r.t.) $\calF_t$, satisfies the following integrability condition:
\begin{align}\label{eq:admissible}
    \EE \int_0^T \theta_t^2 S_t^2 dt < \infty,
\end{align}
the process \eqref{wealth} admits a unique solution, and satisfies $X^{\theta,x_0}_t \geq 0$, $\PP$-almost surely, for any $t \in [0,T]$. We denote by $\Theta(x_0)$ the set of all admissible controls when the initial wealth is $x_0$. 

\subsection{Utility indifference pricing and utility optimization}\label{sec:uti}
Consider a derivative paying $C_T= C(S_T,Y_T)$ at terminal time $T$, for which we can not set a unique price due to the incompleteness of the market. To calculate the indifference price, we define the following value function $V^{(k)}(t,x,s,y)$:
\begin{align*}
    V^{(k)}(t,x,s,y) :=\sup_{\theta\in\Theta(x)}J^{(k)}(t,x,s,y,\theta),
\end{align*}
where
\begin{align}
   J^{(k)}(t,x,s,y,\theta) :=\EE[U(X_T^{\theta}+kC_T)|X_t^{\theta}=x,S_t = s,Y_t = y],
\end{align}
for some function $U:\R\to\R$, called the {\it utility function}. 
This is the maximal utility the investor can achieve with wealth $x$ and $k$ units of the derivative when the value of the tradable asset is $s$ and the value of the auxiliary process is $y$ at time $t$. The function $U$ is the utility function. Then, the indifference price at time $t=0$ for $k$ units of the derivative with $S_0 = s_0$, $Y_0 = y_0$, and initial wealth $x_0$, is defined as the solution $p(k)$ to the equation:
\begin{align}\label{ind-pri}
    V^{(k)}(0,x_0-p(k),s_0,y_0) = V^{(0)}(0,x_0,s_0,y_0).
\end{align}
On the left-hand side we have the utility when the investor buys $k$ holds of the derivative for a price $p(k)$ from her initial wealth $x_0$. On the right-hand side, we have the utility when the investor does nothing (and thus pays nothing from her initial wealth $x_0$). Therefore, to find the indifference price $p(k)$, we need to:
\begin{enumerate}
    \item evaluate $V^{(0)}(0,x_0,s_0,y_0)$;
    \item evaluate the function $V^{(k)}(0,\cdot,s_0,y_0)$;
    \item find the value $p(k)$ such that 
    \begin{align}\label{eq:UID}
    V^{(k)}(0,x_0-p(k),s_0,y_0) = V^{(0)}(0,x_0,s_0,y_0).
    \end{align}
\end{enumerate}

\subsection{The HJB equation for $\phi=V^{(k)}$}\label{sec:hjb}
The following HJB equation for $\phi=V^{(k)}$ is given by:
\begin{align}\label{HJB:og}
\begin{cases}
    \partial_t\phi + H(t,\phi,D\phi,D^2\phi)= 0,\\
    \phi(T,x,s,y) = U(x+kC_T),
    \end{cases}
\end{align}
where 
\begin{align*}
    H(t,\phi,D\phi,D^2\phi) = \sup_{\theta\in\R^m} L^\theta \phi,
\end{align*}
is the Hamiltonian, and where $L^\theta \phi$ is the {\it infinitesimal generator}, which in the setting described here, is given by, 
\begin{align*}
    L^\theta \phi :=(rx+\theta\mu-r\theta s) \phi_x + \frac 12 \phi_{xx}\theta\sig(\theta\sig)^T + bD_y \phi + \frac 12 \text{tr}(a^T (D^2_y \phi) a)+\mu D_s \phi \\+ \frac 12 \text{tr}(\sig^T(D^2_s\phi)\sig) + \theta \sig \rho a^TD^2_{xy} \phi + \theta\sig\sig^TD^2_{xs} \phi+ \sum(\sig\rho a^TD^2_{ys} \phi);
\end{align*}
see, e.g. \cite{survey}. When $k = 0$, the terminal cost is $U(X^\theta_T)$. We refer to such  case as a utility maximization problem.

In the above (and in what follows), we suppressed the variables $(t,s,y)$ from the functions. We use the notation  $\phi_t=\partial_t\phi$, $\phi_x=\partial_x\phi$, and $\phi_{xx}=\partial^2_x\phi$; furthermore, $D_y \phi \in \R^{n\times 1}$ and $D^2_y \phi \in \R^{n\times n}$ stand for the gradient and the Hessian w.r.t.~$y$, $D_s \phi \in \R^{m\times 1}$ and $D^2_s\phi \in \R^{m\times m}$ for the gradient and the Hessian w.r.t.~$s$, and $D^2_{xy} \phi\in \R^{n\times 1}, D^2_{xs}\phi\in \R^{m\times 1}, D^2_{ys}\phi\in \R^{n\times m}$ for the gradients of $\phi_x$ w.r.t.~$y, s$ and the mixed derivatives of $\phi$ w.r.t.~$y$ and $s$.

To evaluate the Hamiltonian, in what follows we assume that 
    the utility function $U(\cdot)$ is concave.

\begin{remark}
    An often-used example of such a function is the negative exponential, i.e., $U(x) = -\gamma^{-1} e^{-\gamma x}$. The concavity of $U$, together with the following convexity property that follows from \eqref{eq:admissible}: whenever  $\theta_1 \in \Theta(x_1)$ and $\theta_2 \in \Theta(x_2)$, then for any $\lambda \in (0,1)$, $\lambda \theta_1 + (1-\lambda)\theta_2 \in \Theta (\lambda x_1 + (1-\lambda) x_2)$, gives the concavity of $V$ in $x$,  one can check Proposition 3.1, \cite{Zari_power}. As a result, $V_{xx} < 0$, the generator becomes quadratic in $\theta$ with a negative quadratic coefficient, so the Hamiltonian can be trivially solved at its critical point,
    and is given by
    \begin{align*}
        H(t,\phi,D\phi,D^2\phi) &= \frac 12 \text{tr}(\sig^T(D^2_s\phi)\sig)+\frac 12 \text{tr}(a^T (D^2_y \phi)a)+\sum(\sig\rho a^T D^2_{ys}\phi)\\
        &\quad+ \mu D_s \phi+rxV_x+bD_y \phi - \frac 12 B^TA^{-1}B,
    \end{align*}
    with
    \begin{align*}
        A = V_{xx}\sig\sig^T  \quad \textrm{and} \quad B = \sig\sig^TD^2_{xs}\phi+\sig \rho a^T D^2_{xy}\phi+(\mu-rs) \phi_x.
    \end{align*}
\end{remark}

Here, we present a conventional verification result, asserting that when a suitably smooth function complies with the HJB equation, it will be equivalent to the value function $V$. For a comprehensive treatment of this result in the general context, we direct the reader to \cite[Theorem 3.5.2]{Huyen}.

\begin{proposition}\label{HJB:thm}
    Consider a function $V^*
    \in \calC^2((0,\infty)^{1+m+n}\times[0,T])$ 
satisfying the terminal condition $$V^*(t,x,s,y) = U(x+kC_T),$$ the quadratic growth condition 
\begin{align*}
    V^*(t,x,s,y) \leq C(1+|x|^2+|s|^2+|y|^2), \quad \forall (t,x,s,y) \in \calC^2((0,\infty)^{1+m+n}\times[0,T]),
\end{align*}
and assume that there is an admissible control $\theta^*$, such that 
    \begin{align*}
        V^*_t + H(t,G,DV^*,D^2V^*) = V^*_t + L^{\theta^*}V^* = 0.
    \end{align*}
    Then
    \begin{align*}
        V^*(t,x,s,y) = J^{(k)}(t,x,s,y,\theta^*) = \sup_{\theta \in \Theta(x)}J^{(k)}(t,x,s,y,\theta)=V^{(k)}(t,x,s,y).
    \end{align*}
  
\end{proposition}

\subsection{Two examples for simplified HJB equations}\label{sec:eg}
We now consider the Black-Scholes model for $S$. Specifically,
\begin{align*}
    dS^i_t &= \mu_i(Y_t)S^i_tdt+S^i_t\sum_j\sigma_{ij}(Y_t)dB^j_t,\quad S^i_0 = s^i_0>0,\quad i = 1,\dots,m,\\
    dY_t &= b(Y_t)dt + a(Y_t)dW_t,\quad Y_0 = y_0 > 0.
\end{align*}
In this case, the value of $S$ actually scales out of the problem, namely, $V^{(k)}(t,x,s,y)$ is constant in $s$. The generator and the Hamiltonian is simplified by removing the derivatives w.r.t. $s$, i.e., we have
\begin{align*}
    L^\theta \phi := &(rx +\theta(\mu \odot s)-r\theta s)\phi_x + \frac 12 \phi_{xx}\theta\sig\odot s(\theta\sig\odot s)^T 
    \\
    &\quad+ bD_y \phi + \frac 12 \text{tr}(a^T (D^2_y \phi) a) + \theta\sig\odot s\rho a^TD^2_{xy} \phi,
\end{align*}
where $\odot$ stands for the component-wise multiplication, i.e.,
\begin{align*}
    \mu \odot s = \begin{bmatrix}
        \mu_1s^1\\\mu_2s^2\\\dots\\\mu_ms^m
    \end{bmatrix}, \quad \sig\odot s =\begin{bmatrix}
         \sig_{11}s^1&\sig_{12}s^1&\dots&\sig_{1m}s^1\\\sig_{21}s^2&\sig_{22}s^2&\dots&\sig_{2m}s^2\\\dots&\dots&\dots&\dots\\\sig_{m1}s^m&\sig_{m2}s^m&\dots&\sig_{mm}s^m
    \end{bmatrix},
\end{align*}
and the HJB equation becomes:
\begin{align*}
    \begin{cases}
    \phi_t+ \frac 12 \text{tr}(a^T (D^2_y \phi) a)+rx\phi_x+bD_y \phi - \frac 12 B^TA^{-1}B=0,\\
    \phi(T,x,y) = U(x+kC_T),
    \end{cases}
    \end{align*}
    with
    \begin{align*}
       A = \phi_{xx}\sig\odot s(\sig\odot s)^T   \quad \textrm{and} \quad B = \sig\odot s\rho a^TD^2_{xy} \phi+(\mu \odot s)-r s) \phi_x.
    \end{align*}
We now provide two examples, where the model could further be specified.

\begin{example}{(Stochastic volatility model)}\label{egSV}

We let the auxiliary process $Y$ denote the stochastic volatility. Consider the dynamics from \cite{invest}:
    \begin{align*}
        &dS_t = (r+Y^1_t)S_tdt + \sig\sqrt{Y_t^2}S_tdB_t,\\
        &\begin{bmatrix}
        dY^1_t\\dY^2_t
    \end{bmatrix} = \begin{bmatrix}\theta^1(k^1-Y^1_t)\\\theta^2(k^2-Y^2_t)
    \end{bmatrix}dt + \begin{bmatrix}1&0\\0&\sqrt{Y^2_t}\end{bmatrix}\begin{bmatrix}a^{11}&a^{12}\\a^{21}&a^{22}
    \end{bmatrix}\begin{bmatrix}dW^1_t\\dW^2_t
    \end{bmatrix},
    \end{align*}
    where $S$ is the risky asset, and $Y^1, Y^2$ are stochastic processes that determine respectively the drift and diffusion coefficient of $S$. In this example, the HJB is given explicitly by
    \begin{align*}
        \begin{cases}
    \phi_t+ \frac 12 \left((a_{11}^2+a_{12}^2)\phi_{y_1y_1}+(2a_{11}a_{21}\sqrt{y_2}+2a_{12}a_{22}\sqrt{y_2})\phi_{y_1y_2}+(a_{21}^2y_2+a_{22}^2y_2)\phi_{y_2y_2}\right)\\\qquad+rx\phi_x+\theta^1(k^1-y_1)\phi_{y_1}+\theta^2(k^2-y_2)\phi_{y_2} \\
    \qquad- \frac 1{2\sig^2y_2\phi_{xx}} \left((\rho_1a_{11}+\rho_2a_{21}\sqrt{y_2})\phi_{xy_1}+(\rho_1a_{12}+\rho_2a_{22}\sqrt{y_2})\phi_{xy_2}+(y_1-r)\phi_x\right)^2=0,\\\phi(T,x,y) = U(x+kC_T).
    \end{cases}
    \end{align*}
\end{example}

\begin{example}{(Non-tradable stocks)}\label{egBS}
We let the auxiliary process $Y$ denote the non-tradable stocks. Consider the dynamics where both $S$ and $Y$ follow constant coefficient Black-Scholes model:
    \begin{align*}
dS_t &= (r+\lambda\sig) S_tdt+\sig S_tdB_t,\quad S_0 = s_0,\\
    dY_t &= b\odot Y_tdt + a\odot Y_tdW_t,\quad Y_0 = y_0 > 0,
\end{align*}
where
\begin{align*}
    \lambda = \frac{\mu-r}{\sig}.
\end{align*}
is the Sharpe ratio for the risky asset $S$, and $\mu,\sig$ are constants. We also have $b$ and $a$ to be constant vector and matrix in $\R^{1\times n}$ and $\R^{n\times n}$, and
\begin{align*}
    b \odot y = \begin{bmatrix}
        b_1y^1\\b_2y^2\\\dots\\b_ny^n
    \end{bmatrix}, \quad a\odot y =\begin{bmatrix}
         a_{11}y^1&a_{12}y^1&\dots&a_{1n}y^1\\a_{21}y^2&a_{22}y^2&\dots&a_{2n}y^2\\\dots&\dots&\dots&\dots\\a_{n1}y^n&a_{n2}y^n&\dots&a_{nn}y^n
    \end{bmatrix}.
\end{align*}
In this example, the HJB is given explicitly by
    \begin{align*}
        \begin{cases}
    \phi_t+ \frac 12 \text{tr}((a\odot y)^T (D^2_y \phi) a\odot y)+rx\phi_x+b\odot yD_y \phi - \frac {(\rho a^TD^2_{xy} \phi+\lambda\phi_x)^2}{2\phi_{xx}} =0,\\
    \phi(T,x,y) = U(x+kC_T).
    \end{cases}
    \end{align*}
\end{example}

\section{Examples}\label{sec:4}

In this section, we pedagogically describe the application of the analytical form of NG, outlined in \ref{sec:2}, to several cases of Example \ref{egBS}. Specifically, In subsection \ref{eg1}, we set $n = 0$, and study the optimal investment problem, demonstrating how a simple choice of nonlinear trial function correctly describes the investment dynamics. In subsection \ref{eg2}, we let $n = 1$, i.e., the case where there is only one non-tradable stock. In subsection \ref{eg3}, we study the case for general $n$. In the example $n \geq 1$, the derivative is a forward, and takes the form $\sum y_i$.  In subsection \ref{cap:figures}, we give explanation and the parameters on the numerical experiments and figures we get.
\subsection{Example 1}\label{eg1}
In order to convert Example \ref{egBS} into an initial value problem of the form \eqref{e:IVP}, we begin by performing the following change of variables,
\begin{equation*}
u(t,x)=-\phi(T-t,x).
\end{equation*}
Since $n=0$, then in terms of the notation introduced in section \ref{sec:2}, the spatial domain $\Omega = [0,\infty)$ is of dimension $d = n+1=1$, and the number of spatial derivatives occurring in the PDE is $\kap=2$. In particular, the explicit form of the right-hand side of the PDE \eqref{e:IVP} is given by, 
\begin{equation*}
    F\big[t,x,u(t,x),u_x(t,x),u_{xx}(t,x)\big] = r x u_x(t,x) - \frac{\lambda^2}{2} \frac{u_x(t,x)^2}{u_{xx}(t,x)},
\end{equation*}
and the initial condition has the following functional form,
\begin{equation*}
    f(x) = 2e^{-\frac{1}{2}x}.
\end{equation*}
Motivated by the form of the initial condition above, we consider an extremely simple trial wave function with a single adjustable parameter. Even though $p=1$ here, we continue to use the boldface notation to represent the parameters $\bm{\theta} \in \mathbb{R}^p$ of the wave function, in order to match with the notation of section \ref{sec:2}. The explicit normalized wave function $\psi_{\bm{\theta}}$ and the associated trial function $u_\theta$ are given as follows,
\begin{align*}
    \psi_{\bm\theta}(x) & := \sqrt{\beta} \, e^{-\frac{1}{2}\beta x}, & x \geq 0, \\
    u_\theta(x) & = \alpha \, \psi_{\bm\theta}(x), & x \geq 0,
\end{align*}
where $\theta_0 = \log \alpha$, $\bm{\theta} = \log \beta$ and $\theta = (\log\alpha,\log \beta) \in \mathbb{R}^2$. Based on the form of the initial condition given above, it is clear that the required choice of initial variational parameters is given by,
\begin{equation*}
    \alpha = 2 \qquad \text{and} \qquad \beta =1.
\end{equation*}
Thus, all that remains to complete the NG for this example is to perform the integrals defining the matrix $M$ and the vector $V$ in \eqref{e:general}, and then to solve the resulting system of ordinary differential equations. The necessary integrals can be performed in closed form and one finds the following results,
\begin{align*}
    M & = \alpha^2
    \begin{bmatrix}
    1 & 0 \\
    0 & \frac{1}{4}
    \end{bmatrix}, \\
    V & = \alpha^2
    \begin{bmatrix}
    -\frac{1}{2}(\lambda^2 + r) \\
    \frac{1}{4}r
    \end{bmatrix}.
\end{align*}
Plugging the above expressions into \eqref{e:general}, we find the following time dependence for the variational parameters,
\begin{equation*}
    \frac{{\rm d}}{{\rm d}t}
    \begin{bmatrix}
    \log\alpha \\
    \log \beta
    \end{bmatrix}
    =
    \begin{bmatrix}
    -\frac{1}{2}(\lambda^2+r) \\
    r
    \end{bmatrix}.
\end{equation*}
Thus, the solution of the ODEs which satisfies the intial conditions is clearly given by
\begin{align*}
    \alpha(t) & = 2e^{-\frac{1}{2}(\lambda^2+r)t}, \\
    \beta(t) & = e^{rt}.
\end{align*}
Remarkably, one finds that the NG reproduces the exact solution of the PDE in this case, which can be verified by plugging $u_{\theta(t)}$ back into \eqref{e:IVP}.

\subsection{Example 2}\label{eg2}
Again, we convert Example \ref{egBS} into an initial value problem of the form \eqref{e:IVP} by performing the following change of variables,
\begin{align*}
    u(t,x,y) = -\phi(T-t,x,y).
\end{align*}
Set $a(y) = a_0y, b(y) = b_0y$ ($a_0=0.3,b_0=0.2,\rho=0.1$).
Then $d = 2, \kap = 5$, and
\begin{align*}
    &F[t,x,y, u, u_x,u_{xx},u_{y},u_{xy},u_{yy}] = \frac12a(y)^2 u_{yy} + rxu_x + b(y) u_{y} -\frac{(\rho a(y) u_{xy} +\lambda u_x)^2}{2u_{xx}},
\end{align*}
\begin{align*}
    f(x,y) = 2e^{-\frac 12 (x+y)}.
\end{align*}
Consider the following trial function with $p=2$ parameters,
\begin{align*}
    u_\theta(x,y) & = \alpha \, \psi_{\bm\theta}(x,y), & x,y \geq 0, \\
    \psi_{\bm\theta}(x,y) & := \sqrt{\beta\zeta} \, e^{-\frac{1}{2}(\beta x+\zeta y)}, & x,y \geq 0,
\end{align*}
where $\theta_0 = \log \alpha$, $\bm{\theta} = (\log \beta,\log \zeta)$ and $\theta = (\log\alpha,\log \beta,\log \zeta) \in \mathbb{R}^3$. The required initial condition for the parameters is thus
\begin{equation*}
    \alpha = 2, \quad \quad \quad \beta =1, \quad \quad \quad \zeta = 1,
\end{equation*}
\begin{align*}
    M
    & =
    \alpha^2
    \begin{bmatrix}
    1 & 0 & 0 \\
    0 & \frac{1}{4} & 0 \\
    0 & 0 & \frac{1}{4}
    \end{bmatrix}, \\
    V
    & =
    \frac{1}{4}
    \alpha^2
    \begin{bmatrix}
    a_0(a_0 + 2\lambda\rho-a_0\rho^2)-2(b_0+r+\lambda^2) \\
    r \\
    b_0-a_0\lambda\rho+a_0^2(\rho^2-1)
    \end{bmatrix}.
\end{align*}
So the parameters evolve as
\begin{align*}
    \alpha(t) & = 2e^{\frac{1}{4}\big[a_0(a_0 + 2\lambda\rho-a_0\rho^2)-2(b_0+r+\lambda^2)\big]t}, \\
    \beta(t) & = e^{rt}, \\
    \zeta(t) & = e^{\big[b_0-a_0\lambda\rho+a_0^2(\rho^2-1)\big]t}.
\end{align*}
Hence, the trial solution is
\begin{equation*}
    u(t,x,y) = \alpha(t)\sqrt{\beta(t)\zeta(t)}e^{-\frac{1}{2}\big[\beta(t)x+\zeta(t)y\big]}.
\end{equation*}
\subsection{Example 3}\label{eg3}
Consider now the $y\in \mathbb{R}^n$ for any $n\geq 1$. In addition, assume that
\begin{align*}
    a(y) & = a_0 \diag(y), \\
    b(y) & = b_0 y.
\end{align*}
As before, we convert to an initial value problem by defining
\begin{align*}
    u(t,x,y) = -\phi(T-t,x,y).
\end{align*}
Then, $d = n+1$, $\kap= \frac{n^2}{2}+\frac 52 n + 2$, and
\begin{align*}
    F =\frac 12 \text{tr}(a^T (D^2_y u)a)+rxu_x +bD_y u - \frac{(\rho a^TD^2_{xy} u+\lambda u_x)^2}{2u_{xx}},
\end{align*}
\begin{align*}
    f(x,y) = 2e^{-\frac 12 (x+\sum y_i)},
\end{align*}
where $D^2_y u$ is the Hessian of $n\times n$, $D_y u$ is the gradient of $n\times 1$, $D^2_{xy} u$ is the gradient of $u_x$ w.r.t. $y$ of $n\times 1$.

Trial function:
\begin{equation*}
    \psi_{\bm{\theta}}(x,y_1,\ldots, y_n) = \sqrt{\beta \zeta^n} e^{-\frac{1}{2}[\beta x + \zeta(y_1+\cdots y_n)]}.
\end{equation*}
Initial condition
\begin{equation*}
    \alpha = 2, \quad \quad \quad \beta =1, \quad \quad \quad \zeta = 1,
\end{equation*}
\begin{align*}
    M
    & =
    \alpha^2\diag(1,1/4,1/4), \\
    V
    & =
    \frac{1}{4}
    \alpha^2
    \begin{bmatrix}
    a_0\left(na_0 + 2n\lambda\rho-\frac{n(n+1)}{2}a_0\rho^2\right)-2(nb_0+r+\lambda^2) \\
    r \\
    n\left[b_0 - a_0\rho\lambda + a_0^2 \left(\frac{n+1}{2}\rho^2 - 1\right) \right]
    \end{bmatrix}.
\end{align*}
So the parameters evolve as
\begin{align}
    \alpha(t) & = 2e^{\frac{1}{4}\big[a_0\left(Na_0 + 2n\lambda\rho-\frac{n(n+1)}{2}a_0\rho^2\right)-2(nb_0+r+\lambda^2)\big]t},   \notag \\
    \beta(t) & = e^{rt} ,\notag \\
    \zeta(t) & = e^{n\left[b_0 - a_0\rho\lambda + a_0^2 \left(\frac{n+1}{2}\rho^2 - 1\right) \right]t}. \label{e:example_3}
\end{align}

\subsection{Discussion of graph}\label{cap:figures}

\begin{figure}
    \centering
    \begin{subfigure}[b]{0.45\textwidth}
          \raggedright
         \includegraphics[width=\textwidth]{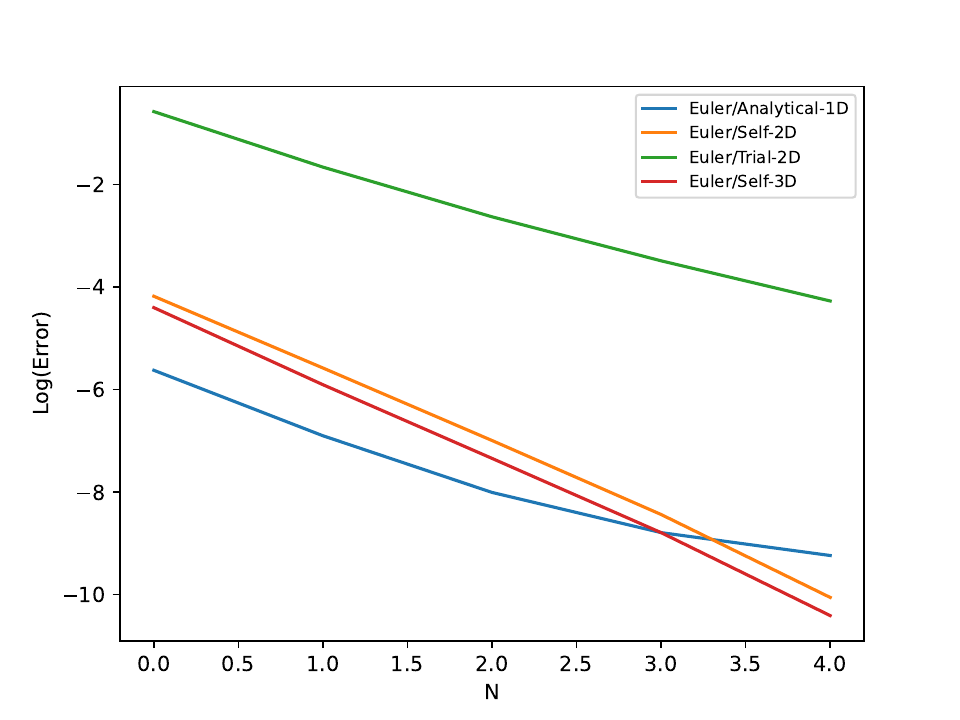}
         \caption{$d = 1, 2, 3$}
    \label{fig:1}
     \end{subfigure}
    \hfill
     \begin{subfigure}[b]{0.45\textwidth}
          \raggedleft
         \includegraphics[width=\textwidth]{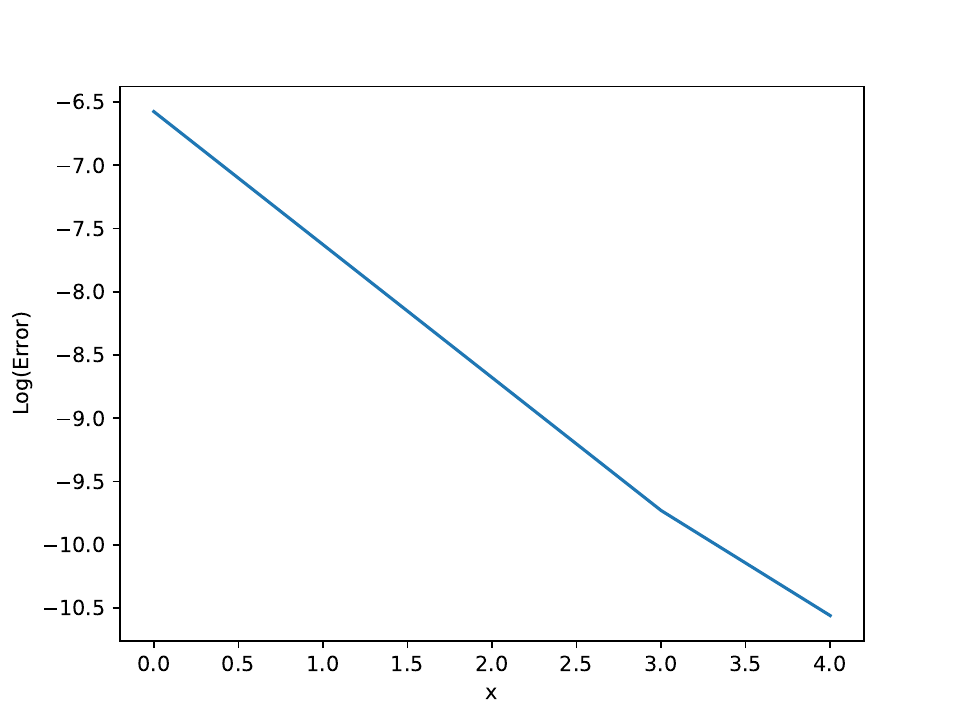}
         \caption{$d = 3$}
    \label{fig:2}
     \end{subfigure}
     \caption{
     \textit{Horizontal axis $N$ is $\log_2$ of the number of space grid points along any dimension (Subfigure (a)) and $x$ is the fixed given $x$ value of the point used to evaluate the error (Subfigure (b)) while the vertical axis is $\log$ absolute error averaged over space points, and $d$ is the the total space dimension including $x, y$.}}

\end{figure}

\begin{figure}
    \centering
    \begin{subfigure}[b]{0.45\textwidth}
         \centering
    \includegraphics[width=\textwidth]{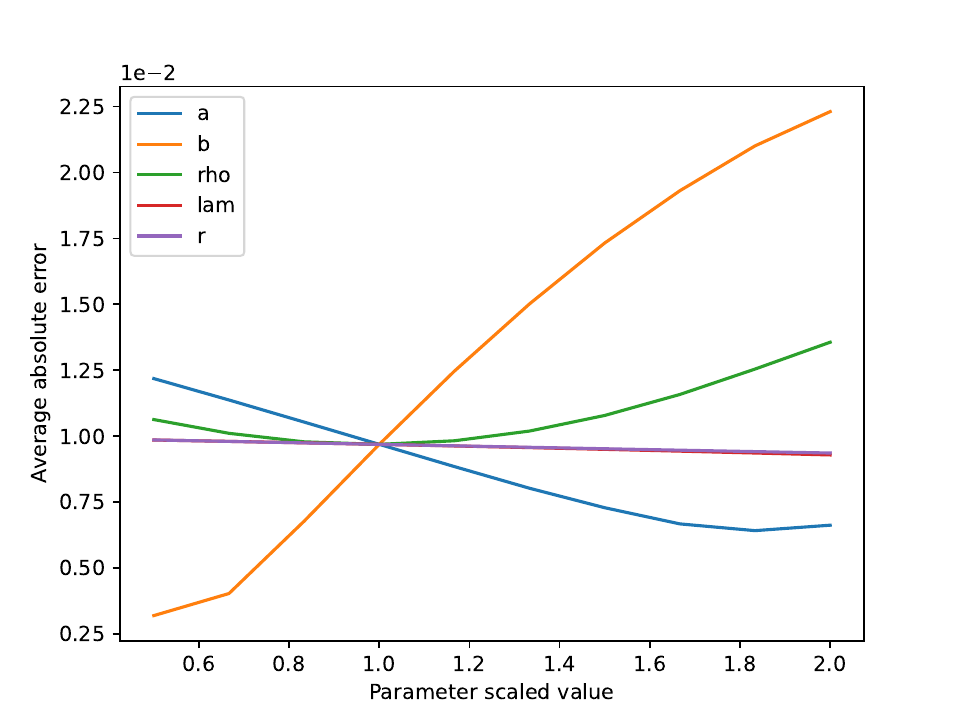}
    \caption{}
    \label{fig:3}
     \end{subfigure}
    \hfill
     \begin{subfigure}[b]{0.45\textwidth}
          \raggedleft
         \includegraphics[width=\textwidth]{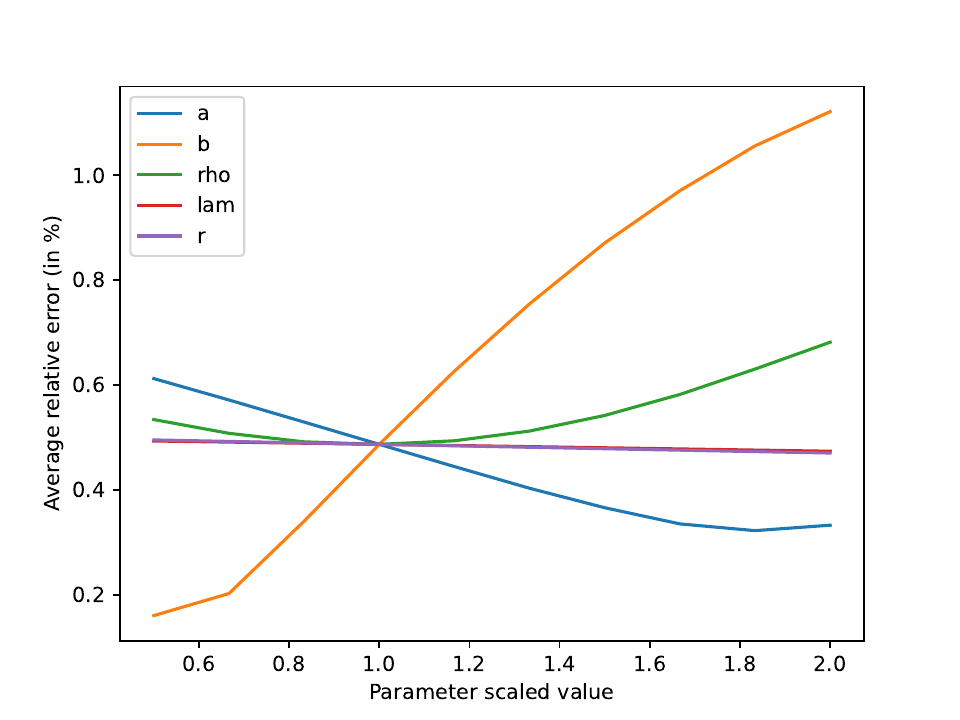}
    \caption{}
    \label{fig:4}
     \end{subfigure}
     \caption{\textit{(a) The absolute error and (b) the relative error averaged over all space points, between forward Euler with $N = 4$ and our trial solution}. 
     }

\end{figure}

We will now present some graphs of numerical results. In Figure \ref{fig:1} and Figure \ref{fig:2}, we aim to show that our trial solution is indeed a good approximation, by computing the error between it and the forward Euler method.

Specifically, for Figure \ref{fig:1}, we have universal parameters, $\gamma = 0.5$, $r = 0.05$ and $\lambda = 0.1$. The errors are all measured at terminal time $T = 1$. The space grid is $[0,4]$ along all dimensions. For the blue curve, the total space dimension $d$ is $1$ (we have no $y$ variable). The error is the $\log$ absolute $\ell_1$ error between forward Euler and analytical solution, averaged over all space grid points. For the orange and green curve, we have additional parameters $a = 0.3$, $b = 0.2$, $\rho = 0.1$. The total space dimension including $x,y$ is $d=2$. The orange curve is the $\log$ absolute $\ell_1$ self-error of forward Euler method, at point $(x,y) = (2,2)$. The green curve is the $\log$ absolute $\ell_2$ error between forward Euler with $N=5$ and trial solution, averaged over all space grid points. For the red curve, we have additional parameters $a =\begin{bmatrix}0.3&0.0\\0.0&0.3\end{bmatrix}$, $b = [0.2,0.2]$, $\rho = [0.1,0.1]$. The total space dimension is $d = 3$. It is the $\log$ absolute $\ell_1$ self-error of forward Euler, at point $(x,y_1,y_2) = (2,2,2)$.

For Figure \ref{fig:2}, we have the parameters $\gamma = 0.5$, $r = 0.05$, $\lambda = 0.1$, $a =\begin{bmatrix}0.3&0.0\\0.0&0.3\end{bmatrix}$, $b = [0.2,0.2]$ and $\rho = [0.1,0.1]$. The error is measured at terminal time $T = 1$. The total space dimension including $x,y$ is $d=3$, and the space grid is $[0,4]$ along all dimensions. On the vertical axis, we have the $\log$ absolute $\ell_2$ error between forward Euler with $N=4$ and our trial solution, averaged over all $(y_1,y_2)$ grids, but at fixed given $x$ value, specified as the horizontal axis.

From these figures, we have the following observations. The blue, orange and red curves in Figure \ref{fig:1} indicate that when $d = 1,2,3$, the forward Euler method will converge, as we increase the number of grid points along each space dimension. The green curve in Figure \ref{fig:1} indicates that the forward Euler method also has decreasing errors with our trial function, as we increase the number of grid points. This shows that our trial solution is within a certain range of the true solution, which is unknown. In Figure \ref{fig:2}, we are comparing when $d = 3$, how the error between our trial solution and forward Euler varies according to $x$. We can see that the error tends to decrease as $x$ is approaching upper bound, which in terms of the context, is as the price of stock is increasing.

Now, in Figures \ref{fig:3} and \ref{fig:4}, we study how changing the values of parameters effects the errors. On the vertical axis, we have the average percentage $\ell_1$ relative error and average $\ell_1$ absolute error, averaged over all space grid points. The total space dimension including $x,y$ is $d=3$, and the space grid is $[0,4]$ along all dimensions. For different curves, we vary one parameter at each time. The default value of parameters are $T = 1$, $\gamma = 0.5$, $r = 0.05$, $\lambda = 0.1$, $a =\begin{bmatrix}0.3&0.0\\0.0&0.3\end{bmatrix}$, $b = [0.2,0.2]$ and $\rho = [0.1,0.1]$. On the horizontal axis, we use a scale of $[0.5,2]$ to indicate the relative value of each parameter in their respective range. Write  $a =\begin{bmatrix}a_0&0.0\\0.0&a_0\end{bmatrix}$, $b = [b_0,b_0]$ and $\rho = [\rho_0,\rho_0]$, the ranges for each parameter are $a_0:[0.25,0.4]$, $b_0:[0.1,0.4]$, $\rho_0: [-0.5,0.4]$, $\lambda: [0.05,0.2]$, $r: [0.25,0.1]$. Specifically, when the horizontal axis takes value $1$, it means the parameter takes the default value. 

We see from the graphs that varying $\lambda$ and $r$ has little influence on the error, changing $a$ and $\rho$ affects the error a bit, while the error grows dramatically as increasing the value of $b$. 
 
\section{Conclusions and future directions}\label{sec:conclusion}
In summary, we have introduced a family of nonlinearly parametrized trial functions which provide approximate analytical solutions to the HJB equation and serve as a precursor to the Neural Galerkin scheme. Although the trial functions introduced here may appear simple, they effectively capture the time-dependence of the solutions, achieving exact precision in certain cases. In cases where analytically tractable trial functions only capture the overarching qualitative behavior, it is plausible that minor deformations, when parameterized by neural networks, could elevate them to high precision. Future work will focus on learning deformations induced by parametrized diffeomorphisms of the spatial coordinate, treating the analytical trial function as the base distribution of a normalizing flow-based neural network \cite{papamakarios2021normalizing}.
\section{Acknowledgements} Authors gratefully acknowledge support from NSF under grants DMS-2038030 and DMS-2006305.
\section{Technical Proofs}\label{sec:5}
\begin{proof}[Proof of \eqref{e:example_3}]
Plugging the trial function 
\begin{equation*}
    u_\theta = \alpha \sqrt{ab^n} e^{-\frac{1}{2}[a x + b (y_1 + \cdots y_n)]}.
\end{equation*}
into the right-hand side of the PDE we obtain the following expression (in a shorthand notation),
\begin{align*}
    F_\theta 
    & = \left[\frac{1}{8} a_0^2 b^2 (y_1^2 + \cdots y_n^2) - \frac{1}{2} r a x - \frac{1}{2} b_0 b (y_1 + \cdots + y_n) 
    - \frac{\big[\frac{1}{2}\rho a_0 a b (y_1+\cdots+y_n) - \lambda a\big]^2}{2a^2}\right] u_\theta \notag \\
    & = \left[\frac{1}{8} a_0^2 b^2 \sum_{i,j}(\delta_{ij}-\rho^2)y_i y_j + \frac{1}{2}(\lambda\rho a_0 - b_0) b\sum_{i} y_i - \frac{1}{2}rax - \frac{\lambda^2}{2}\right] u_\theta.
\end{align*}
Expressing the trial function in terms of the variational parameters,
\begin{equation*}
    u_\theta = e^{\theta_0 + \frac{1}{2}(\theta_1 + n\theta_2)-\frac{1}{2}[e^{\theta_1} x + e^{\theta_2} (y_1 + \cdots y_n)]},
\end{equation*}
we obtain the following variational derivatives,
\begin{align*}
    \frac{\partial u_\theta}{\partial\theta_0} & = u_\theta, \\
    \frac{\partial u_\theta}{\partial\theta_1} & = \frac{1}{2}(1 - a x) u_\theta, \\
    \frac{\partial u_\theta}{\partial\theta_2} & = \frac{1}{2}[n - b (y_1+\cdots+y_n)] u_\theta.
\end{align*}
It will be useful to keep in mind the following identities,
\begin{align}
    \langle \psi_{\bm\theta} | x |\psi_{\bm\theta} \rangle & = \frac{1}{a} \label{e:identity_x}, \\
    \langle \psi_{\bm\theta} | y_i |\psi_{\bm\theta} \rangle & = \frac{1}{b} \label{e:identity_x}, \\
    \langle \psi_{\bm\theta} | y_i y_j |\psi_{\bm\theta} \rangle & = \frac{\delta_{ij}+1}{b^2} \label{e:identity_yy}, \\
    \langle \psi_{\bm\theta}| \sum_{i,j}y_i^2 y_j | \psi_{\bm\theta} \rangle \label{e:identity_y2y}
    & = \frac{2n(n+2)}{b^3}, \\
    \langle \psi_{\bm\theta}| \sum_{i,j,k}y_i y_j y_k | \psi_{\bm\theta} \rangle \label{e:identity_yyy}
    & = \frac{(n+2)(n+1)n}{b^3}.
\end{align}
Thus,
\begin{align*}
    \frac{1}{\alpha^2} \langle u_\theta | F_\theta \rangle
    & = \frac{1}{8} a_0^2 b^2 \sum_{i,j}(\delta_{ij}-\rho^2) \langle \psi_{\bm\theta} | y_i y_j |\psi_{\bm\theta}\rangle + \frac{1}{2}(\lambda\rho a_0 - b_0) b\sum_{i} \langle \psi_{\bm\theta} | y_i | \psi_{\bm\theta} \rangle 
    \\&\quad- \frac{1}{2}ra\langle \psi_{\bm\theta} |x | \psi_{\bm\theta} \rangle - \frac{\lambda^2}{2} \langle \psi_{\bm\theta} | \psi_{\bm\theta} \rangle  \notag \\
    & = \frac{1}{8} a_0^2 b^2 \sum_{i,j}(\delta_{ij}-\rho^2)(\delta_{ij}+1) + \frac{n(\lambda\rho a_0 - b_0) - r - \lambda^2}{2} \\
    & = \frac{a_0^2 n [2-\rho^2(n+1)]}{8} + \frac{n(\lambda\rho a_0 - b_0) - r - \lambda^2}{2},
\end{align*}
\begin{align*}
    \frac{1}{\alpha^2}\left\langle \frac{\partial u_\theta}{\partial \theta_0} \middle| F_\theta \right\rangle
    & = \frac{1}{\alpha^2} \langle u_\theta | F_\theta \rangle, \\
    \frac{1}{\alpha^2}\left\langle \frac{\partial u_\theta}{\partial \theta_1} \middle| F_\theta \right\rangle
    & = \frac{1}{2\alpha^2}\left[\langle u_\theta | F_\theta \rangle - a \langle u_\theta | x | u_\theta \rangle \right] \\
    & = -\frac{ra}{4}\left[\langle \psi_{\bm\theta} | x | \psi_{\bm\theta} \rangle - a \langle \psi_{\bm\theta} | x^2 | \psi_{\bm\theta} \rangle \right] \\
    & = \frac{r}{4}, \\
    \frac{1}{\alpha^2}\left\langle \frac{\partial u_\theta}{\partial \theta_2} \middle| F_\theta \right\rangle
    & = \frac{1}{2\alpha^2}\left[n \langle u_\theta | F_\theta \rangle - b \langle u_\theta | \sum_i y_i | u_\theta \rangle \right] \\
    & = \frac{1}{2}\bigg[n\left(\frac{a_0^2 n [2-\rho^2(n+1)]}{8} + \frac{n(\lambda\rho a_0 - b_0)}{2}\right) + \notag \\
    & \quad - b\bigg(\frac{1}{8} a_0^2 b^2 \sum_{i,j,k}(\delta_{ij}-\rho^2) \langle \psi_{\bm\theta} | y_i y_jy_k |\psi_{\bm\theta}\rangle + \frac{1}{2}(\lambda\rho a_0 - b_0) b\sum_{i,j} \langle \psi_{\bm\theta} | y_i y_j | \psi_{\bm\theta} \rangle \bigg) \bigg] \\
    & = \frac{1}{2}\bigg[n\left(\frac{a_0^2 n [2-\rho^2(n+1)]}{8} + \frac{n(\lambda\rho a_0 - b_0)}{2}\right) + \notag \\
    & \quad - b\bigg\{\frac{1}{8} a_0^2 b^2 \bigg(\frac{2n(n+2)}{b^3} - \rho^2 \frac{(n+2)(n+1)n}{b^3} \bigg) + \frac{n(n+1)(\lambda\rho a_0 - b_0)}{2b} \bigg\}\bigg] \\
    & = \frac{n}{4}\left[b_0 - a_0\rho\lambda + a_0^2 \left(\frac{n+1}{2}\rho^2 - 1\right) \right].
\end{align*}
\end{proof}



\bibliographystyle{plain}
\bibliography{refs}
\end{document}